\documentclass[12pt]{amsart}

\usepackage{amsmath,amsfonts,amsthm}

\usepackage{fullpage}
\usepackage{times,txfonts}
\newtheorem{theorem}{Theorem}
\newtheorem{prop}[theorem]{Proposition}
\newtheorem{coro}[theorem]{Corollary}
\newtheorem{lemma}[theorem]{Lemma}

\DeclareMathOperator{\sh}{sh}
\DeclareMathOperator{\ch}{ch}

\begin{document}

\renewcommand{\thefootnote}{\fnsymbol{footnote}}

\baselineskip 17pt
\parskip 7pt
\sloppy




\title[Half-Derivatives of the Andrews--Gordon Identity]{${q}$-Series
      and
      ${L}$-Functions
      Related to Half-Derivatives of
      the      Andrews--Gordon      Identity
    }


    \author{Kazuhiro \textsc{Hikami}}


  \address{Department of Physics, Graduate School of Science,
    University of Tokyo,
    Hongo 7--3--1, Bunkyo, Tokyo 113--0033, Japan.
    }

    \email{\texttt{hikami@phys.s.u-tokyo.ac.jp}}


\date{November 21, 2002, Revised on May 8, 2003}

\begin{abstract}
Studied is  a generalization of Zagier's $q$-series identity.
We introduce a generating function of $L$-functions at non-positive
integers, which is regarded as a half-differential of the
Andrews--Gordon $q$-series.
When $q$ is a root of unity, 
the generating function coincides with the quantum
invariant for the torus knot.
\end{abstract}





\maketitle

\section{Introduction}

In Ref.~\cite{DZagie01a}, Zagier studied
the $q$-series,
\begin{equation}
  X(q) =
  \sum_{n=0}^\infty
  (1-q) \, (1-q^2) \cdots (1-q^n)
\end{equation}
and proved that the asymptotic expansion is given by
\begin{equation}
  X(\mathrm{e}^{-t})
  =
  \mathrm{e}^{t/24}
  \sum_{n=0}^\infty 
  \frac{T(n)}{n!} \left( \frac{t}{24} \right)^n
  \label{Zagier_identity}
\end{equation}
Here $T(n)$ is the Glaisher $T$-number,
\begin{equation*}
  \frac{\sh(2 \,  x)}{\ch(3 \, x)}
  =
  \sum_{n=0}^\infty
  \chi_{12}(n) \, \mathrm{e}^{-nx}
  =
  2 \sum_{n=0}^\infty
  (-1)^n
  \,
  \frac{T(n)}{(2 \, n+1)!} \, x^{2n+1}
\end{equation*}
and is given in terms of the Dirichlet $L$-function as
\begin{equation}
  T(n)
  =
  \frac{1}{2} \, (-1)^{n+1} \, L(-2 \, n -1, \chi_{12})
\end{equation}
where $\chi_{12}(n)$ is the Dirichlet character with modulus $12$
defined by
\begin{equation*}
  \begin{array}{c|ccccc}
    n \mod 12 & 1 & 5 & 7 & 11 & \text{others}
    \\
    \hline
    \chi_{12}(n) & 1 & -1 & -1 & 1 & 0
  \end{array}
\end{equation*}
It was pointed out that the right hand side of
eq.~\eqref{Zagier_identity} is regarded as a half-differential of the
Dedekind $\eta$-function with weight $1/2$.
Interesting is that the function $X(q)$ is intimately connected with
the knot theory;
it is a generating function of an upper bound of the number of
linearly independent
Vassiliev invariants.

Purpose of this paper is to 
study a generalization of Zagier's identity
(see  Refs.~\cite{AndrUrroOnok01a,CoogKOno02a,LoveKOno02a} for
this attempt).
Our motivation
is based on an observation that
the $q$-series $X(q)$ with $q$ being  root of unity
appears as a colored Jones invariant
of the trefoil~\cite{Kasha95,MuraMura99a}.
We shall show that the $q$-series, which reduces to the invariant of
the torus knot in a case of $q$ being root of unity,
becomes the generating function of the $L$-function with negative
integers.
We note that
a relationship between the modular form and the quantum invariant
was  discussed in Ref.~\cite{LawrZagi99a}, where the
Witten--Reshetikhin--Turaev invariant of the Poincar\'e homology
sphere was studied.
Throughout this paper
we use a standard notation,
\begin{gather*}
  (x)_n
  = (x;q)_n
  =
  \prod_{i=1}^n
  (1- x \, q^{i-1}) 
  \\[2mm]
  (x_1,\dots,x_j)_n
  =
  (x_1,\dots,x_j;q)_n
  =
  (x_1)_n \cdots (x_j)_n
  \\[2mm]
  \begin{bmatrix}
    n \\
    c
  \end{bmatrix}
  =
  \begin{cases}
    \displaystyle
    \frac{(q)_n}{(q)_c \, (q)_{n-c}}  
    & \text{for $0 \leq c \leq n$}
    \\[3mm]
    0 & \text{otherwise}
  \end{cases}
\end{gather*}

We state the main result of this article.
Let the $q$-series $X_2^{(a)}(q)$ for $a=0,1$ be
\begin{align}
  X_2^{(0)}(q)
  & =
  \sum_{n=0}^\infty \sum_{c=0}^n
  (q)_n  \, q^{c^2+c} \,
  \begin{bmatrix}
    n \\
    c
  \end{bmatrix}
  \\[2mm]
  X_2^{(1)}(q)
  & =
  \sum_{n=1}^\infty \sum_{c=0}^n
  (q)_{n-1}  \, q^{c^2} \,
  \begin{bmatrix}
    n \\
    c
  \end{bmatrix}
\end{align}

\begin{theorem}
  \label{thm:RR}
  \begin{align}
    X_2^{(0)}(\mathrm{e}^{-t})
    & =
        \mathrm{e}^{9 t/40}
    \sum_{n=0}^\infty
    \frac{T_2^{(0)}(n)}{n!} \,
    \left(
      \frac{t}{40}
    \right)^n
    \label{X20_e}
    \\[2mm]
    X_2^{(1)}(\mathrm{e}^{-t})
    & =
        \mathrm{e}^{ t/40}
    \sum_{n=0}^\infty
    \frac{{T}_2^{(1)}(n)}{n!} \,
    \left(
      \frac{t}{40}
    \right)^n
    \label{X21_e}
  \end{align}
  where
  \begin{equation*}
    T_2^{(a)}(n)
    =
    \frac{1}{2} \, (-1)^{n+1} \, L(-2 \, n -1, \chi_{20}^{(a)})
  \end{equation*}
\end{theorem}
We note that, using the Mellin transformation, we have the generating
function of the $T$-series;
\begin{equation*}
  \frac{\sh(2 \, (a+1) \, x)}{\ch(5 \, x)}
  =
  \sum_{n=0}^\infty
  \chi_{20}^{(a)}(n) \, \mathrm{e}^{-nx}
  =
  2 \sum_{n=0}^\infty
  (-1)^n  \,
  \frac{T_2^{(a)}(n)}{(2 \, n+1)!} \,
  x^{2n+1}
\end{equation*}
where the periodic function $\chi_{20}^{(a)}(n)$ are
\begin{gather*}
  \begin{array}{c|ccccc}
    n \mod 20  
    &    3    &     7    &     13    &    17
    &\text{other}
    \\
    \hline
    {\chi}_{20}^{(0)}(n) &    1 &   -1 & -1 &    1 & 0
  \end{array}
  \\[2mm]
  \begin{array}{c|ccccc}
    n \mod 20  
    &    1    &     9    &     11    &    19
    &\text{other}
    \\
    \hline
    {\chi}_{20}^{(1)}(n) &
    1 &   -1 & -1 &    1 & 0
  \end{array}
\end{gather*}

To see a relationship with a modular form,
we recall the well known
Rogers--Ramanujan identity (see, \emph{e.g.},
Refs.~\cite{Andre76,Andre86,GaspRahm90});
\begin{subequations}
  \begin{align}
    (q)_\infty \,
    \sum_{n=0}^\infty
    \frac{q^{n^2+n}}{(q)_n}
    & =
    (q, q^4, q^5 ; q^5)_\infty
    \nonumber \\
    & =
    \sum_{n=0}^\infty
    \chi_{20}^{(0)}(n) \, q^{(n^2-9)/40}
    \label{RR}
    \displaybreak[0]
    \\[2mm]
    (q)_\infty \,
    \sum_{n=0}^\infty
    \frac{q^{n^2}}{(q)_n}
    & =
    (q^2, q^3, q^5 ; q^5)_\infty
    \nonumber \\
    & =
    \sum_{n=0}^\infty
    {\chi}_{20}^{(1)}(n) \, q^{(n^2-1)/40}
    \label{RR2}
  \end{align}
\end{subequations}
Those functions are the two-dimensional representation of the modular
group with weight $1/2$,
and Theorem~\ref{thm:RR} indicates that the $q$-series
$X_2^{(a)}(q)$ is related with
a  ``half-differential'' of the Rogers--Ramanujan $q$-series;
\begin{equation}
  \label{half_RR}
  X_2^{(a)}(q)
  =
  - \frac{1}{2} \,
  \sum_{n=0}^\infty
  n \, \chi_{20}^{(a)}(n) \, q^{\frac{n^2 - (3-2a)^2}{40}}
\end{equation}
Here two sides cannot be defined simultaneously but the equality
holds as the Taylor expansions of $q\to \mathrm{e}^{-t}$.

The Rogers--Ramanujan identity can be generalized to the
Andrews--Gordon identity~\cite{Andre86};
let $m \in \mathbb{Z}_{>1}$ and $0 \leq a \leq m-1$, then
\begin{align}
  \sum_{n_1 \geq \dots \geq n_{m-1} \geq 0}^\infty
  \frac{
    q^{n_1^{~2} + n_2^{~2} + \dots + n_{m-1}^{~2}
      + n_{a+1} + \dots + n_{m-1}}
  }{
    (q)_{n_1 - n_2} (q)_{n_2 - n_3} \dots (q)_{n_{m-1}}
  }
  & =
  \prod_{\substack{
      n=1 \\
      n \not\equiv 0, \pm (a+1) \mod 2m+1
    }}^\infty
  (1-q^n)^{-1}
  \nonumber
  \\
  & =
  \frac{1}{(q)_\infty} \,
  \sum_{n=0}^\infty
  \chi_{8m+4}^{(a)}(n) \,
  q^{
    \frac{n^2 - (2m-2a-1)^2}{8(2m+1)}
  }
\end{align}
where
we have introduced the periodic function $\chi_{8m+4}^{(a)}(n)$ as
\begin{equation*}
  \begin{array}{c|ccccc}
    n \mod (8 \, m+4)
    & 2 \, m - 2 \, a -1  & 2 \, m + 2 \, a +3
    & 6 \, m -2 \, a + 1  & 6 \, m + 2 \, a + 5 & \text{others}
    \\
    \hline
    \chi_{8m+4}^{(a)}(n) & 1 & -1 & -1 & 1 & 0
  \end{array}
\end{equation*}
As a generalization of Theorem~\ref{thm:RR}
to the $q$-series related to a half-derivative of the
Andrews--Gordon $q$-series in a sense of eq.~\eqref{half_RR},
we define
the function $X_m^{(a)}(q)$ with
$m \in \mathbb{Z}_{>0}$ and $a=0,1,\dots,m-1$ by
\begin{equation}
  X_m^{(a)}(q)
  =
  \sum_{k_1,k_2,\dots,k_m=0}^\infty
  (q)_{k_m} \,
  q^{k_1^{~2} + \dots + k_{m-1}^{~2} + k_{a+1} + \dots + k_{m-1}}
  \left(
    \prod_{\substack{i=1 \\ i \neq a}}^{m-1}
    \begin{bmatrix}
      k_{i+1} \\
      k_i
    \end{bmatrix}
  \right) \cdot
  \begin{bmatrix}
    k_{a+1} + 1 \\
    k_a
  \end{bmatrix}
  \label{define_X_a}
\end{equation}

\begin{theorem}
  \label{thm:general_m}
  \begin{equation}
    X_m^{(a)}(\mathrm{e}^{-t})
    =
    \mathrm{e}^{
      \frac{(2m -2a-1)^2}{8 (2m+1)} t
    }
    \sum_{n=0}^\infty
    \frac{T_m^{(a)}(n)}{n!} \,
    \left(
      \frac{t}{8 \, (2\,m+1)}
    \right)^n
  \end{equation}
  where $T$-series is given by the $L$-function
  \begin{align}
    \label{T_and_L_function}
    T_m^{(a)}(n)
    & =
    \frac{1}{2} \, (-1)^{n+1} \,
    L(-2 \, n -1 , \chi_{8m+4}^{(a)})
    \\
    & =
    (-1)^n \, 2^{4n} \,
    \frac{(2 \, m+1)^{2n+1}}{n+1} \,
    \sum_{r=1}^{8m+4}
    \chi_{8m+4}^{(a)}(r) \,
    B_{2n+2}
    \left(
      \frac{r}{8 \, m+4}
    \right)
    \nonumber
  \end{align}
  with the Bernoulli polynomial $B_n(x)$.
\end{theorem}
In this case the generating function of the $T$-series  is written as
\begin{equation}
  \frac{\sh \bigl( 2 \, (a+1) \, x \bigr)}
  {\ch \bigl( (2\, m+1) \, x \bigr)}
  =
  \sum_{n=0}^\infty
  \chi_{8m+4}^{(a)}(n) \, \mathrm{e}^{-nx}
  =
  2 \sum_{n=0}^\infty
  (-1)^n \,
  \frac{T_m^{(a)}(n)}{(2 \, n+1)!} \, x^{2n+1}
\end{equation}
See that
a case of $m=2$ corresponds to Theorem~\ref{thm:RR} and
that  $m=1$ is nothing but
Zagier's identity~\eqref{Zagier_identity}.
Furthermore above Theorem shows that
\begin{equation}
  \label{half_AG}
  X_m^{(a)}(q)
  =
  - \frac{1}{2}
  \sum_{n=0}^\infty
  n \, \chi_{8m+4}^{(a)}(n) \,
  q^{
    \frac{n^2 - (2m-2a-1)^2}{8(2m+1)}
  }
\end{equation}
as a generalization of eq.~\eqref{half_RR}.

For our later convention  we collect $q$-series identities as follows
(see, \emph{e.g.}, Refs.~\cite{Andre76,GaspRahm90});
\begin{itemize}
\item 
  $q$-binomial coefficient
  \begin{subequations}
    \begin{align}
      \begin{bmatrix}
        n+1 \\
        c
      \end{bmatrix}
      & =
      q^c \,
      \begin{bmatrix}
        n \\
        c
      \end{bmatrix}
      +
      \begin{bmatrix}
        n \\
        c-1
      \end{bmatrix}
      \label{binomial_1}
      \\
      & =
      \begin{bmatrix}
        n \\
        c
      \end{bmatrix}
      +
      q^{n+1-c} \,
      \begin{bmatrix}
        n \\
        c-1
      \end{bmatrix}
      \label{binomial_2}
    \end{align}
  \end{subequations}
  
\item
  $q$-binomial formula
  \begin{equation}
    \label{binomial_3}
    \sum_{n=0}^\infty
    \frac{(a)_n}{(q)_n} \, z^n
    =
    \frac{(a \, z)_\infty}{( z)_\infty}
  \end{equation}

\item
  the Jacobi triple product identity
  \begin{equation}
    \label{Jacobi}
    \sum_{k= - \infty}^\infty
    (-1)^k q^{\frac{1}{2} k^2}  x^k
    =
    (q)_\infty \, 
    (x^{-1} q^{\frac{1}{2}})_\infty \,
    (x \,  q^{\frac{1}{2}})_\infty
  \end{equation}
\end{itemize}


In Section~\ref{sec:proof1} we prove Theorem~\ref{thm:RR}.
Section~\ref{sec:proof2} is for the proof of
Theorem~\ref{thm:general_m}.
Strategy to prove these theorems is essentially same with a proof of
eq.~\eqref{Zagier_identity} in Ref.~\cite{DZagie01a};
we define the function $H_m^{(a)}(x)$ and derive the $q$-series
$X_m^{(a)}(q)$ as a differential of $H_m^{(a)}(x)$.
We comment on a relationship between the quantum knot invariant and
our $q$-series in Section~\ref{sec:knot}.

\section*{Acknowledgement}
The author would like to thank Hitoshi Murakami for bringing
Ref.~\cite{DZagie01a}
to  attention.
He  thanks Ken Ono for useful communications.
Thanks are also to George Andrews for sending Ref.~\cite{GAndre01a}.
He thanks Don Zagier for suggesting a proof of Theorem~\ref{conj:main}.
This work is supported in part by the Sumitomo foundation.

\section{Proof of Theorem~\ref{thm:RR}}
\label{sec:proof1}

We define
\begin{align}
  \label{define_H}
  H(x)
  & =
  \sum_{n=0}^\infty \sum_{c=0}^n
  (x)_{n+1} \, x^n \cdot
  q^{c^2+c} \, x^{2 c}
  \begin{bmatrix}
    n \\
    c
  \end{bmatrix}
\end{align}

\begin{prop}
  \label{prop:RR1}
  \begin{equation}
    \label{RR_x}
    H(x)
    =
    \sum_{n=0}^\infty
    \chi_{20}^{(0)}(n) \,
    q^{(n^2-9)/40} \,
    x^{(n-3)/2}
  \end{equation}
\end{prop}

\begin{proof}
  We prove this identity by showing that both hand sides satisfy the
  same difference equation.
  For the right hand side, we have
  \begin{align}
    H(x)
    & =
    \sum_{n=0}^\infty
    \chi_{20}(n) \,
    q^{(n^2-9)/40} \,
    x^{(n-3)/2}
    \nonumber
    \\
    & =
    1 - q \, x^2 +
    \sum_{n=10}^\infty
    \chi_{20}(n) \,
    q^{(n^2-9)/40} \, x^{(n-3)/2}
    \nonumber
    \\
    & =
    1 - q \, x^2 - q^4 \, x^5 \, H(q\,x)
    \label{diff_eq_H}
  \end{align}
  where we have used $\chi_{20}(n+10) = - \chi_{20}(n)$.

  To study the difference equation for the left hand side,
  we further define
  \begin{equation}
    H(x,y,z)
    =
    \sum_{n=0}^\infty \sum_{c=0}^n
    (x)_{n} \, y^n \cdot
    q^{c^2+c} \, z^{2 c}
    \begin{bmatrix}
      n \\
      c
    \end{bmatrix}
  \end{equation}
  and we investigate the difference equation of $H(x,y,z)$.


  We have
  \begin{align}
    H(x,y,q\, z)
    & =
    \sum_{n=0}^\infty \sum_{a=1}^{n+1}
    (x)_n \, y^n \, q^{a^2 + a-2} \, z^{2 a -2} \,
    \begin{bmatrix}
      n \\
      a-1
    \end{bmatrix}
    \displaybreak[0]
    \nonumber
    \\
    &
    \stackrel{\eqref{binomial_1}}{=}
    \sum_{m=1}^\infty \sum_{a=0}^m
    (x)_{m-1} \, y^{m-1} \, q^{a^2+a-2} \, z^{2a-2}
    \begin{bmatrix}
      m \\
      a
    \end{bmatrix}
    \nonumber
    \\
    & \qquad
    - q^{-2} z^{-2}
    \sum_{n=0}^\infty \sum_{a=0}^n
    (x)_n \, y^n \, q^{a^2+a} \,
    (q^{1/2} z )^{2 a}
    \begin{bmatrix}
      n \\
      a
    \end{bmatrix}
    \displaybreak[0]
    \nonumber
    \\
    & =
    \sum_{m=1}^\infty \sum_{a=0}^m
    \frac{(q^{-1} x)_{m}}{1-q^{-1} x} \,
    y^{m-1} \, q^{a^2+a-2} \, z^{2a-2} \,
    \begin{bmatrix}
      m \\
      a
    \end{bmatrix}
    \nonumber \\
    & \qquad \qquad \qquad
    - q^{-2} z^{-2} \, H(x, y, q^{1/2} z)
    \displaybreak[0]
    \nonumber 
    \\
    & =
    \frac{y^{-1} q^{-2} z^{-2}}{1- q^{-1} x} \,
    \left(
      H(q^{-1} x , y, z) - 1
    \right)
    - q^{-2} z^{-2} \, H(x, y, q^{1/2} z)      
    \label{difference_2}
  \end{align}

  In the same manner, we have following;
  \begin{align}
    H(x, q\, y, z)
    & =
    \sum_{n=0}^\infty \sum_{c=0}^n
    (x)_n \,
    \bigl( 1- (1-q^n \,x ) \bigr) \, x^{-1} \,
    y^n \, q^{c^2+c} \, z^{2c}
    \begin{bmatrix}
      n \\
      c
    \end{bmatrix}
    \displaybreak[0]
    \nonumber
    \\
    & \stackrel{\eqref{binomial_2}}{=}
    x^{-1} H(x,y,z)
    -
    \sum_{n=0}^\infty \sum_{c=0}^{n+1}
    (x)_{n+1} \,  x^{-1} \, y^n \,  q^{c^2+c} \,  z^{2c}
    \begin{bmatrix}
      n+1 \\
      c
    \end{bmatrix}
    \nonumber
    \\
    & \qquad
    +
    \sum_{n=0}^\infty \sum_{c=1}^{n+1}
    (x)_{n+1} \,  x^{-1} \, y^n \,  q^{c^2+n+1} \, z^{2c}
    \begin{bmatrix}
      n \\
      c-1
    \end{bmatrix}
    \displaybreak[0]
    \nonumber
    \\
    & =
    x^{-1} H(x,y,z)
    -
    \sum_{m=1}^\infty \sum_{c=0}^{m}
    (x)_{m} \,  x^{-1} \, y^{m-1} \, q^{c^2+c} \, z^{2c}
    \begin{bmatrix}
      m \\
      c
    \end{bmatrix}
    \nonumber
    \\
    & \qquad
    +
    \sum_{n=0}^\infty \sum_{a=0}^{n}
    (x)_{n+1} \,  x^{-1} \, y^n \,  q^{a^2+2a+n+2} \, z^{2a+2}
    \begin{bmatrix}
      n \\
      a
    \end{bmatrix}
    \displaybreak[0]
    \nonumber
    \\
    & =
    x^{-1} \,  y^{-1} + x^{-1} \,  (1-y^{-1}) \, H(x,y,z)
    +
    q^2 \, (1-x) \, x^{-1} \, z^2
    \, H(q \, x , q\, y, q^{1/2} z)
    \label{difference_3}
  \end{align}

  We combine these two difference equations;
  from eqs.~\eqref{difference_3}
  we eliminate $H(q \, x, q\, y,q^{\frac{1}{2}} z)$
  using eq.~\eqref{difference_2}.
  We get
  \begin{multline}
    (x \, y -q \, z^2) \, H(x, q\, y , z)
    + ( 1- y) \, H(x,y,z)
    \\
    =
    - q^4 \, (1-x) \, y \, z^4 \, H(q \, x , q\, y, q \, z)
    +
    1- q \, z^2
    \label{diff_medium}
  \end{multline}
  When we substitute
$
    (x,y,z) \to (q \, x , x, x)
$ to this equation, the first term vanishes.
  Recalling that
  \begin{equation}
    H(x)
    = (1-x) \, H(q \,x , x, x)
  \end{equation}
  we obtain eq.~\eqref{diff_eq_H}.
\end{proof}

Setting $x\to 1$ in Prop.~\ref{prop:RR1}, we see that the right hand
side reduces to the Rogers--Ramanujan $q$-series
due to the Jacobi triple identity,
\begin{equation*}
  H(x=1)=(q,q^4,q^5;q^5)_\infty
\end{equation*}
For the left hand side~\eqref{define_H}, we see this fact
from the
following lemma,
which can be proved by use of the binomial formula~\eqref{binomial_3}.
\begin{lemma}
  \begin{equation}
    H(x)
    =
    (q \, x)_\infty \,
    \sum_{c=0}^\infty \frac{q^{c^2+c}}{(q \, x)_c} \,
    x^{3c}
    +
    (1-x)
    \sum_{n=0}^\infty \sum_{c=0}^n
    \bigl(
    (q \, x)_n - (q \, x)_\infty
    \bigr) \, x^n
    \, q^{c^2+c} \, x^{2c} \,
    \begin{bmatrix}
      n \\
      c
    \end{bmatrix}
    \label{RR_2}
  \end{equation}
\end{lemma}

Before proceeding to the proof of Theorem~\ref{thm:RR}, we recall
the known  result on the Mellin transformation;
\begin{prop}
  \label{prop:Mellin}
  Let $\chi_p$ be a periodic function with modulus $p$ with mean value
  zero,
  and
  \begin{equation*}
    L(s,\chi_p)
    =
    \sum_{n=1}^\infty \frac{\chi_p(n)}{n^s}
  \end{equation*}
  As $t \searrow 0$ we have
  \begin{equation}
    \sum_{n=0}^\infty
    n \, \chi_p(n) \,  \mathrm{e}^{- n^2 t}
    \sim
    \sum_{n=0}^\infty L(-2 \, n -1 , \chi_p) \, \frac{(-t)^n}{n!}
  \end{equation}
\end{prop}

\begin{proof}
  Assumptions of $\chi_p$ support that $L(s,\chi_p)$ has an analytic
  continuation to $\mathbb{C}$.
  We apply the Mellin transformation to
  \begin{equation*}
    \sum_{n=0}^\infty
    n \, \chi_p(n) \,  \mathrm{e}^{- n^2 t}
    \sim
    \sum_{n=0}^\infty \gamma_n \, t^n
  \end{equation*}
  From the left hand side, we have
  \begin{align*}
    \sum_{n=0}^\infty n \, \chi_p(n) \,
    \int_0^\infty t^{s-1} \, \mathrm{e}^{-n^2 t} \,
    \mathrm{d} t
    & =
    \sum_{n=0}^\infty  \chi_p(n) \,
    \frac{\Gamma(s)}{n^{2s-1}}
    \\
    & = \Gamma(s) \, L(2 \, s-1, \chi_p)
  \end{align*}
  We also have from the right hand side that
  \begin{align*}
    \int_0^\infty
    \left(
      \sum_{n=0}^{N-1} \gamma_n \, t^n + \mathcal{O}(t^N)
    \right) \,
    t^{s-1} \, \mathrm{d} t
    & =
    \sum_{n=0}^{N-1} \frac{\gamma_n}{n+s} + R_N(s)
  \end{align*}
  where $R_N(s)$ is analytic in $\Re(s)>-N$.
  Thus $\gamma_n$ is the residue of $\Gamma(s) \, L(2 \, s-1,\chi_p)$
  at $s=-N$, and we get
  \begin{equation*}
    \gamma_n =
    (-1)^n \frac{L(-2 \, n -1, \chi_p)}{n!}
  \end{equation*}
\end{proof}


\begin{proof}[Proof of eq.~\eqref{X20_e} in Theorem~\ref{thm:RR}]
We equate eq.~\eqref{RR_x} with eq.~\eqref{RR_2}, and
we set $x\to 1$ after differentiating with respect to $x$.
Using eq.~\eqref{RR}, we get
\begin{multline}
  \label{dH_for_RR1}
  \sum_{n=0}^\infty
  \frac{n}{2} \,
  \chi_{20}^{(0)}(n) \, q^{(n^2-9)/40}
  =
  -
  \sum_{n=0}^\infty \sum_{c=0}^n
  \bigl( (q)_n - (q)_\infty \bigr) \, q^{c^2+c} \,
  \begin{bmatrix}
    n \\
    c
  \end{bmatrix}
  \\
  +
  (q)_\infty \sum_{c=0}^\infty
  \frac{q^{c^2+c}}{(q)_c} \,
  \left(
    3 \, c +
    \sum_{i=1}^c
    \frac{q^i}{1-q^i}
  \right)
  +
  (q, q^4, q^5; q^5)_\infty
  \left( \frac{3}{2} -
    \sum_{i=1}^\infty
    \frac{q^i}{1-q^i}
  \right)
\end{multline}

We substitute $q= \mathrm{e}^{-t}$ to eq.~\eqref{dH_for_RR1}, and
study the Taylor expansion of $t$.
Therein terms including
infinite products such as
$(q)_\infty$ and $(q,q^4,q^5;q^5)_\infty$ vanish as they induce an
infinite order of $t$.
Then we get eq.~\eqref{half_RR}.
Using
Prop.~\ref{prop:Mellin} we recover
eq.~\eqref{X20_e} in Theorem~\ref{thm:RR}.
\end{proof}

For the proof of the rest of Theorem~\ref{thm:RR}, we define the
function $G(x)$ by
\begin{equation}
  \label{define_G}
  G(x)
  =
  \sum_{n=1}^\infty \sum_{c=0}^n
  (x)_{n} \, x^{n-1} \cdot
  q^{c^2} \, x^{2 c}
  \begin{bmatrix}
    n \\
    c
  \end{bmatrix}
\end{equation}

\begin{prop}
  \begin{equation}
    \label{RR_y}
    G(x)
    =
    \sum_{n=0}^\infty
    {\chi}_{20}^{(1)}(n) \,
    q^{(n^2-1)/40} \,
    x^{(n-1)/2}
  \end{equation}
\end{prop}

\begin{proof}
  We show that both hand sides satisfy the same $q$-difference equation.
  It is easy to see that the right hand side satisfies
  \begin{equation}
    \label{dif_eq_G_x}
    G(x) = 1- q^2 \, x^4 - q^3 \, x^5 \, G(q \, x)
  \end{equation}
  For the left hand side, we substitute
  $
  (x,y,z) \to (x, x, q^{-1/2} x)
  $
  in eq.~\eqref{diff_medium}.
  Recalling that the function $G(x)$ in eq.~\eqref{define_G} is given
  by
  \begin{equation*}
    G(x) = \frac{1}{x}  \,
    \left(
      H(x, x, q^{-1/2} \, x) - 1
    \right)
  \end{equation*}
  we obtain eq.~\eqref{dif_eq_G_x}.
\end{proof}

One sees from eq.~\eqref{RR_y} using eq.~\eqref{RR2} that
the function $G(x)$ gives the Rogers--Ramanujan $q$-series
\begin{equation*}
  G(x=1) =      (q^2, q^3, q^5 ; q^5)_\infty
\end{equation*}
For expression~\eqref{define_G},
we can rewrite as follows using the binomial
formula~\eqref{binomial_3}.
This lemma supports above equality.
\begin{lemma}
  \begin{multline}
    \label{G_x}
    G(x)
    =
    (q \, x)_\infty \,
    \sum_{c=0}^\infty
    \frac{q^{c^2}}{(q \, x)_c} \, x^{3 c -1}
    -(1-x) \, x^{-1} \, (q \, x)_\infty
    \\
    +
    (1-x) \,
    \sum_{n=1}^\infty \sum_{c=0}^n
    \bigl(
    (q \, x)_{n-1} - (q \, x)_\infty
    \bigr) \, 
    x^{n-1} \, q^{c^2} \, x^{2 c}
    \,
    \begin{bmatrix}
      n \\
      c
    \end{bmatrix}
  \end{multline}
\end{lemma}

\begin{proof}[Proof of eq.~\eqref{X21_e} in Theorem~\ref{thm:RR}]
  We differentiate both eqs.~\eqref{RR_y} and~\eqref{G_x} w.r.t. $x$ and
  substitute $x\to 1$.
  Using eq.~\eqref{RR2} we have
  \begin{multline}
    \sum_{n=0}^\infty
    \frac{n}{2} \,  {\chi}_{20}^{(1)}(n) \, q^{(n^2-1)/40}
    =
    -\sum_{n=1}^\infty \sum_{c=0}^n
    \bigl(
    (q)_{n-1} - (q)_\infty
    \bigr) \,
    q^{c^2} \,
    \begin{bmatrix}
      n \\
      c
    \end{bmatrix}
    \\
    - (q^2, q^3, q^5; q^5)_\infty \,
    \left(
      \frac{1}{2} + \sum_{i=1}^\infty \frac{q^i}{1-q^i}
    \right)
    + (q)_\infty
    + (q)_\infty \sum_{c=0}^\infty
    \frac{q^{c^2}}{(q)_c} 
    \left(
      3 \, c+ \sum_{i=1}^c \frac{q^i}{1-q^i}
    \right)
  \end{multline}
  We substitute $q=\mathrm{e}^{-t}$, and we obtain eq.~\eqref{half_RR}.
  Prop.~\ref{prop:Mellin} proves
  eq.~\eqref{X21_e} in Theorem~\ref{thm:RR}.
\end{proof}

\section{Proof of Theorem~\ref{thm:general_m}}
\label{sec:proof2}

We define the function $H_m^{(a)}(x)$ for $m \in \mathbb{Z}_{>0}$ and
$a=0,1,\dots,m-1$ by
\begin{multline}
  H_m^{(a)}(x)
  =
  \sum_{k_1, \dots, k_m=0}^\infty
  (x)_{k_m+1} \, x^{k_m}
  \,
  \left(
    \prod_{i=1}^{a-1} q^{k_i^{~2}} \, x^{2 k_i} \,
    \begin{bmatrix}
      k_{i+1} \\
      k_i
    \end{bmatrix}
  \right)
  \\
  \times
  q^{k_a^{~2}} \, x^{2 k_a} \,
  \begin{bmatrix}
    k_{a+1}+1 \\
    k_a
  \end{bmatrix}
  \cdot
  \left(
    \prod_{i=a+1}^{m-1}
    q^{k_i^{~2}+k_i} \,
    x^{2 k_i} \,
    \begin{bmatrix}
      k_{i+1} \\
      k_i
    \end{bmatrix}
  \right)
  \label{define_H_m_a}
\end{multline}

\begin{prop}
  \begin{equation}
    \label{H_m_a_and}
    H_m^{(a)}(x)=
    \sum_{n=0}^\infty
    \chi_{8m+4}^{(a)}(n) \,
    q^{
      \frac{n^2 - (2m-2a-1)^2}{8(2m+1)}
    } \,
    x^{\frac{n-(2m-2a-1)}{2}}
  \end{equation}
\end{prop}

\begin{proof}
  Method is essentially same with a proof of Prop.~\ref{prop:RR1};
  we prove that both sides satisfy the same difference equation as a
  function of $x$.
  Anti-periodicity,
  $\chi_{8m+4}^{(a)}(n+4 \, m+2)
  = -  \chi_{8m+4}^{(a)}(n)$,
  shows that
  the r.h.s. satisfies the difference equation
  \begin{equation}
    H_m^{(a)}(x)
    =
    1- q^{a+1} \, x^{2 a +2}
    -
    q^{2m -a} \, x^{2m+1} \,
    H_m^{(a)}( q \, x)
    \label{difference_general_a}
  \end{equation}

  For the l.h.s. we prepare several difference equations
  for the following functions;
  \begin{multline}
    H_m^{(a)}(x, y, z_1, \dots, z_{m-1})
    =
    \sum_{k_1, \dots, k_m=0}^\infty
    (x)_{k_m} \, y^{k_m}
    \\
    \times
    \left(
      \prod_{i=1}^{a-1}
      q^{k_i^{~2}} z_i^{~ 2 k_i} \,
      \begin{bmatrix}
        k_{i+1} \\
        k_i
      \end{bmatrix}
    \right) \cdot
    q^{k_a^{~2}} z_a^{~ 2 k_a} \,
    \begin{bmatrix}
      k_{a+1}+1 \\
      k_a
    \end{bmatrix}
    \\
    \times
    \left(
      \prod_{i=a+1}^{m-1}
      q^{k_i^{~2}+k_i} z_i^{~ 2 k_i} \,
      \begin{bmatrix}
        k_{i+1} \\
        k_i
      \end{bmatrix}
    \right)
  \end{multline}
  We note that
  \begin{equation}
    \label{H_xyz_H}
    H_m^{(a)}(x)
    =
    (1 -x) \, H_m^{(a)}(q \, x, x, \boldsymbol{x})
  \end{equation}
  where for brevity we have used a notation,
  $\boldsymbol{x}=(\underbrace{x,\dots,x}_{m-1})$.

  By applying eq.~\eqref{binomial_1}
  to
  $
  \begin{bmatrix}
    k_{a+1}+1 \\
    k_a
  \end{bmatrix}
  $ in  the definition~\eqref{define_H_m_a}
  of $H_m^{(a)}(x,y,z_1,\dots,z_{m-1})$,
  we obtain  a  following  equation;
  \begin{multline}
    H_m^{(a)}(x,y,z_1,\dots,z_{m-1})
    =
    H_m^{(0)}(x,y, q^{-\frac{1}{2}} z_1, \dots, q^{-\frac{1}{2}}
    z_{a-1},
    z_a,\dots,z_{m-1})
    \\
    +
    q\, z_a^{~2} \,
    H_m^{(a-1)}(x,y,z_1,\dots,z_{a-1},q^{\frac{1}{2}} z_a,
    z_{a+1} ,\dots,z_{m-1})
    \label{gen_dif_1}
  \end{multline}
  When we apply eq.~\eqref{binomial_2} in place of
  eq.~\eqref{binomial_1},
  we get for $a=0,\dots,m-2$
  \begin{multline}
    H_m^{(a)}(x,y,z_1,\dots,z_{m-1})
    =
    H_m^{(0)}(x,y, q^{-\frac{1}{2}} z_1, \dots,
    q^{-\frac{1}{2}}   z_{a},
    z_{a+1},\dots,z_{m-1})
    \\
    +
    q\, z_a^{~2} \,
    H_m^{(a-1)}(x,y,z_1,\dots,z_{a},q^{\frac{1}{2}} z_{a+1},
    z_{a+2} ,\dots,z_{m-1})
    \label{gen_dif_2}
  \end{multline}
  For $a=m-1$ we have
  \begin{multline}
    H_m^{(m-1)}(x,y,
    z_1 , \dots ,
    z_{m-1})
    \\
    =
    H_m^{(0)}(x,y,
    q^{-\frac{1}{2}} z_1,\dots, q^{-\frac{1}{2}} z_{m-1})
    +
    q \, z_{m-1}^{~2} \,
    H_m^{(m-2)}(x, q\, y ,  z_1 , \dots, z_{m-1})
    \label{gen_dif_4}
  \end{multline}

  Next we have
  \begin{align}
    & H_m^{(m-1)}(q \, x, y, z_1, \dots,z_{m-1})
    \nonumber \\
    &=
    \sum_{k_m=1}^\infty
    \sum_{k_{m-1}=0}^{k_m}
    \dots
    \sum_{k_1=0}^{k_2}
    (q \, x)_{k_m-1} \, y^{k_m-1} \,
    \left(
      \prod_{i=1}^{m-1} q^{k_i^{~2}} \, z_i^{~2 k_i} \,
      \begin{bmatrix}
        k_{i+1} \\
        k_i
      \end{bmatrix}
    \right)
    \nonumber \\
    & =
    \frac{1}{(1-x) \, y}
    \left(
      H_m^{(0)} (x, y, q^{-\frac{1}{2}} z_1, \dots ,
      q^{-\frac{1}{2}} z_{m-1})
      -
      1
    \right)
    \label{gen_dif_3}
  \end{align}

  Further we have
  \begin{align}
    & H_m^{(0)}(q \, x , y , z_1, \dots, z_{m-1})
    \nonumber \\
    & =
    \sum_{k_m=0}^\infty
    \sum_{k_{m-1}=0}^{k_m} \dots
    \sum_{k_{1}=0}^{k_2} 
    (q \, x)_{k_m-1} \, (1 - q^{k_m} x) \, y^{k_m} \,
    \left(
      \prod_{i=1}^{m-1} q^{k_i^{~2} + k_i} z_i^{~2 k_i} \,
      \begin{bmatrix}
        k_{i+1} \\
        k_i
      \end{bmatrix}
    \right)
    \nonumber 
    \\
    & =
    \frac{1}{1-x} \,
    \left(
      H_m^{(0)}(x,y,z_1,\dots,z_{m-1})
      - x \,
      H_m^{(0)}(x, q\, y, z_1 , \dots, z_{m-1})
    \right)
    \label{gen_dif_5}
  \end{align}
  We use these difference equations to prove
  Theorem~\ref{thm:general_m}.
  
  A recursive use of eq.~\eqref{gen_dif_1} gives
  \begin{multline}
    H_m^{(m-1)}(q \, x , y, q^{\frac{1}{2}} z_1,\dots,
    q^{\frac{1}{2}}z_{m-1})
    =
    H_m^{(0)}(q \, x, y, z_1,\dots,z_{m-2},q^{\frac{1}{2}} z_{m-1})
    \\
    +
    q^2 \, z_{m-1}^{~2} \,
    H_m^{(0)}(q \, x, y, z_1,\dots,z_{m-3},q^{\frac{1}{2}} z_{m-2},
    q  \, z_{m-1})
    \\
    +
    q^4 \, z_{m-2}^{~2 } \, z_{m-1}^{~2} \,
    H_m^{(0)}(q \, x, y, z_1,\dots,z_{m-4},q^{\frac{1}{2}} z_{m-3},
    q  \, z_{m-2}, q\, z_{m-1})
    \\
    + \dots
    +
    q^{2m-2} \, z_1^{~2} \cdots z_{m-1}^{~2} \,
    H_m^{(0)}(q \, x , y, q \, z_1, \dots, q \, z_{m-1})
    \label{gen_dif_1_1}
  \end{multline}
  while a recursive use of eq.~\eqref{gen_dif_2} gives
  \begin{multline}
    H_m^{(m-2)}(q \, x , y, q^{\frac{1}{2}} z_1,\dots,
    q^{\frac{1}{2}}z_{m-1})
    =
    H_m^{(0)}(q \, x, y, z_1,\dots,z_{m-2},q^{\frac{1}{2}} z_{m-1})
    \\
    +
    q^2 \, z_{m-2}^{~2} \,
    H_m^{(0)}(q \, x, y, z_1,\dots,z_{m-3},q^{\frac{1}{2}} z_{m-2},
    q  \, z_{m-1})
    \\
    +
    q^4 \, z_{m-3}^{~2 } \, z_{m-2}^{~2} \,
    H_m^{(0)}(q \, x, y, z_1,\dots,z_{m-4},q^{\frac{1}{2}} z_{m-3},
    q  \, z_{m-2}, q\, z_{m-1})
    \\
    + \dots
    +
    q^{2m-4} \, z_1^{~2} \cdots z_{m-2}^{~2} \,
    H_m^{(0)}(q \, x , y,  q^{\frac{1}{2}} z_1,
    q \, z_2, \dots, q \, z_{m-1})
    \label{gen_dif_2_1}
  \end{multline}

  Then we get
  \begin{align*}
    & H_m^{(0)}(x,y,z_1,\dots,z_{m-1})
    - x \, H_m^{(0)}(x, q\, y, z_1,\dots,z_{m-1})
    \\
    &
    \stackrel{\eqref{gen_dif_5}}{=}
    (1-x) \, H_m^{(0)}(q \, x, y, z_1, \dots ,z_{m-1})
    \displaybreak[0]
    \\
    &
    \stackrel{\eqref{gen_dif_4}}{=}
    (1-x) \,
    \Bigl(
      H_m^{(m-1)}(q \, x, y, q^{\frac{1}{2}} z_1,\dots,
      q^{\frac{1}{2}}z_{m-1})
      \\
    &  \qquad \qquad \qquad
      -
      q^2 \, z_{m-1}^{~2} \,
      H_m^{(m-2)}(q \, x, q\, y, q^{\frac{1}{2}}z_1,
      \dots, q^{\frac{1}{2}} z_{m-1})
    \Bigr)
    \displaybreak[0]
    \\
    &
    \stackrel{\eqref{gen_dif_3}}{=}
    y^{-1} \left(
      H_m^{(0)}(x,y,z_1,\dots,z_{m-1}) -1
    \right)
    \\
    & \qquad \qquad \qquad
    -q^2 (1-x) \, z_{m-1}^{~2} \,
    H_m^{(m-2)}(q \,x , q\,y, q^{\frac{1}{2}} z_1,
    \dots , q^{\frac{1}{2}} z_{m-1})
  \end{align*}
  Substituting eq.~\eqref{gen_dif_2_1} into above equation, we get
  a difference equation for $H_m^{(0)}(x,y,z_1,\dots,z_{m-1})$;
  \begin{multline}
    (1- y^{-1}) \, H_m^{(0)}(x, y, z_1,\dots, z_{m-1})
    -x \, H_m^{(0)}(x , q \, y, z_1, \dots, z_{m-1})
    + y^{-1}
    \\
    =
    -q^2 \, (1-x) \, z_{m-1}^{~2} \,
    \Bigl(
    H_m^{(0)}(q \, x , q\, y, z_1, \dots, z_{m-2} ,
    q^{\frac{1}{2}} z_{m-1})
    \\
    + q^2  \, z_{m-2}^{~2} \,
    H_m^{(0)}(q \, x , q\, y, z_1, \dots, z_{m-3} ,
    q^{\frac{1}{2}} z_{m-2}, q \, z_{m-1})
    \\
    + \dots
    +
    q^{2m-4} \, z_1^{~2} \cdots z_{m-2}^{~2} \,
    H_m^{(0)}(q \, x , q \, y, q^{\frac{1}{2}} z_1,
    q \, z_2, \dots, q \, z_{m-1})
    \Bigr)
    \label{gen_k_8}
  \end{multline}

  Therewith
  by substituting eq.~\eqref{gen_dif_1_1} into eq.~\eqref{gen_dif_3}, we
  obtain another difference equation for
  $H_m^{(0)}(x,y,z_1,\dots,z_{m-1})$;
  \begin{multline}
    \frac{1}{(1-x) \, y} \,
    \left( H_m^{(0)}(x,y,z_1, \dots, z_{m-1})
      -1
    \right)
    -
    q^{2m-2} z_1^{~2} \cdots z_{m-1}^{~2} \,
    H_m^{(0)}(q \, x, y, q \, z_1, \dots, q \, z_{m-1})
    \\
    =
    H_m^{(0)}(q \, x ,  y, z_1, \dots, z_{m-2} ,
    q^{\frac{1}{2}} z_{m-1})
    + q^2  \, z_{m-1}^{~2} \,
    H_m^{(0)}(q \, x ,  y, z_1, \dots, z_{m-3} ,
    q^{\frac{1}{2}} z_{m-2}, q \, z_{m-1})
    \\
    + \dots
    +
    q^{2m-4} \, z_2^{~2} \cdots z_{m-1}^{~2} \,
    H_m^{(0)}(q \, x ,  y, q^{\frac{1}{2}} z_1,
    q \, z_2, \dots, q \, z_{m-1})
    \label{gen_k_9}
  \end{multline}

  We set
  $z_1=z_2=\dots=z_{m-1}=z$ in eqs.~\eqref{gen_k_8}
  and~\eqref{gen_k_9}.
  We can eliminate the right hand side of eq.~\eqref{gen_k_9}
  by use of eq.~\eqref{gen_k_8}, and we
  get
  \begin{multline}
    \label{gen_dif_6}
    -(1-y) \, H_m^{(0)}(x, y, \boldsymbol{z}) + 1 - q \, z^2
    \\
    =
    (x \, y - q \, z^2) \, H_m^{(0)}(x, q\, y, \boldsymbol{z})
    +
    q^{2m} \, (1 -x) \, y \, z^{2m} \,
    H_m^{(0)}(q \, x , q \, y , q \, \boldsymbol{z})
  \end{multline}
  Setting
  $(x,y,z) \to (q\, x, x, x)$, and recalling eq.~\eqref{H_xyz_H},
  we find $H_m^{(0)}(x)$ satisfies $q$-difference
  equation~\eqref{difference_general_a} with $a=0$.


  In the case of  $a \neq 0$, we first recall that
  \begin{multline}
    H_m^{(a)}(x , y, z_1,\dots,z_{m-1})
    =
    H_m^{(0)}(x, y, q^{-\frac{1}{2}} z_1, \dots, q^{-\frac{1}{2}}
    z_{a-1},
    z_a, \dots, z_{m-1})
    \\
    +
    q \, z_a^{~2} \,
    H_m^{(0)}(x, y, q^{-\frac{1}{2}}z_1, \dots, q^{-\frac{1}{2}}
    z_{a-2} , z_{a-1} , q^{\frac{1}{2}} z_a , z_{a+1} , \dots,
    z_{m-1})
    \\
    + \dots
    +
    q^a \, z_1^{~2} \cdots z_a^{~2} \,
    H_m^{(0)}(x,y,
    q^{\frac{1}{2}} z_1,\dots, q^{\frac{1}{2}} z_a,
    z_{a+1}, \dots , z_{m-1})
    \label{H_0_and_a}
  \end{multline}
  which is given by an iterated use of eq.~\eqref{gen_dif_1}.
  We rewrite this identity as
  \begin{equation}
    H_m^{(a)}(x,y,z_1,\dots,z_{m-1})
    =
    \Bigl(
    \Hat{\mathcal{D}}_m^{(a)} \, H_m^{(0)}
    \Bigr) \,
    (x,y,z_1,\dots,z_{m-1})
  \end{equation}
  Here the difference operator is defined by
  \begin{multline}
    \Hat{\mathcal{D}}_m^{(a)}
    =
    \Hat{T}_1^{~-1} \cdots \Hat{T}_{a-1}^{~-1}
    + q \, z_a^{~2} \, \Hat{T}_1^{~-1} \cdots \Hat{T}_{a-2}^{~-1} \,
    \Hat{T}_{a}
    + q^2 \, z_{a-1}^{~2} \, z_a^{~2} \,
    \Hat{T}_1^{~-1} \cdots \Hat{T}_{a-3}^{~-1} \,
    \Hat{T}_{a-1} \,  \Hat{T}_{a}
    \\
    + \dots
    +
    q^a \, z_1^{~2} \cdots z_a^{~2} \, \Hat{T}_1 \cdots \Hat{T}_a
  \end{multline}
  where we have used the $q$-shift operator
  \begin{equation}
    \bigl(
    \Hat{T}_k^{~\pm 1} \, f
    \bigr) \,
    (z_1, \dots, z_{m-1})
    =
    f(z_1, \dots, q^{\pm \frac{1}{2}} z_k , \dots, z_{m-1})
  \end{equation}
  It can be seen by a direct computation that
  for $1 \leq b \leq a$ we have
  \begin{multline}
    \Hat{\mathcal{D}}_m^{(a)} \,
    z_b^{~2} \cdots z_a^{~2} \,
    H_m^{(0)}(x,y, z_1,\dots,z_{b-1},
    q^{\frac{1}{2}} z_b, q \, z_{b+1}, \dots , q\,z_a,
    z_{a+1} , \dots, z_{m-1})
    \\
    =
    q^{b-a} z_b^{~2} \cdots z_a^{~2} \,
    H_m^{(a)}(x, y, z_1,\dots,z_{b-1},
    q^{\frac{1}{2}} z_b, q \, z_{b+1} , \dots , q \, z_a,
    q_{a+1}, \dots , q_{m-1})    
    \label{gen_k_11}
  \end{multline}

  We apply $\Hat{\mathcal{D}}_m^{(a)}$ to eq.~\eqref{gen_k_8}.
  Using
  eq.~\eqref{gen_k_11}, we get
  \begin{multline}
    (1-y^{-1}) \, H_m^{(a)}(x,y,z_1,\dots,z_{m-1})
    -
    x \, H_m^{(a)}(x, q\, y, z_1, \dots, z_{m-1})
    \\
    + y^{-1}(1+q \, z_a^{~2} + \dots + q^a \, z_1^{~2} \dots z_a^{~2})
    \\
    =
    - q^2\,(1-x) \, z_{m-1}^{~2} \,
    \Bigl(
    H_m^{(a)}(q \, x, q\, y, z_1,\dots,z_{m-2},q^{\frac{1}{2}}z_{m-1})
    + \dots +
    \\
    +
    q^{2(m-a-1)} \,  z_{a}^{~2} \cdots z_{m-2}^{~2} \,
    H_m^{(a)}(q \,x , q\, y, z_1, \dots, z_{a-1} ,
    q^{\frac{1}{2}} z_a , q \, z_{a+1}, \dots, q \, z_{m-1})
    \\
    +
    q^{2(m-a)-1} \, z_{a-1}^{~2} \cdots z_{m-2}^{~2} \,
    H_m^{(a)}(q \,x , q\, y, z_1, \dots, z_{a-2} ,
    q^{\frac{1}{2}} z_{a-1} , q \, z_{a}, \dots, q \, z_{m-1})
    \\
    + \dots +
    q^{2m-a-1} \,  z_1^{~2} \cdots z_{m-2}^{~2} \,
    H_m^{(a)}(q \, x, q\, y, q^{\frac{1}{2}} z_1, q \, z_2, \dots,
    q \, z_{m-1})
    \Bigr)
    \label{gen_na_1}
  \end{multline}

  Using eq.~\eqref{H_0_and_a}, eq.~\eqref{gen_k_9} can be rewritten as
  \begin{multline*}
    \frac{1}{(1-x) \, y} \,
    \left( H_m^{(0)}(x,y,z_1, \dots, z_{m-1})
      -1
    \right)
    \\
    =
    H_m^{(0)}(q \, x ,  y, z_1, \dots, z_{m-2} ,
    q^{\frac{1}{2}} z_{m-1})
    +
    q^2  \, z_{m-1}^{~2} \,
    H_m^{(0)}(q \, x ,  y, z_1, \dots, z_{m-3} ,
    q^{\frac{1}{2}} z_{m-2}, q \, z_{m-1})
    \\
    +
    \dots
    +
    q^{2(m-a-2)} \, z_{a+2}^{~2} \cdots z_{m-1}^{~2} \,
    H_m^{(0)}(q \, x ,  y,  z_1, \dots, z_a,
    q^{\frac{1}{2}} z_{a+1},
    q \, z_{a+2}, \dots, q \, z_{m-1})
    +
    \\
    q^{2(m-a-1)} \, z_{a+1}^{~2} \cdots z_{m-1}^{~2} \,
    H_m^{(a)}(q \, x ,  y, q^{\frac{1}{2}} z_1, \dots,
    q^{\frac{1}{2}} z_a,
    q \, z_{a+1}, \dots, q \, z_{m-1})
  \end{multline*}
  Applying $\Hat{\mathcal{D}}_m^{(a)}$ to above equation, we get
  \begin{multline}
    \frac{1}{(1-x) \, y} \,
    \left( H_m^{(a)}(x,y,z_1, \dots, z_{m-1})
      -(1+ q \, z_a^{~2}
      +
      \dots
      + q^a \, z_1^{~2} \cdots z_a^{~2})
    \right)
    \\
    =
    H_m^{(a)}(q \, x, y, z_1, \dots, z_{m-2}, q^{\frac{1}{2}} z_{m-1})
    +
    q^2 z_{m-1}^{~2} \,
    H_m^{(a)}(q \, x, y, z_1, \dots, z_{m-3} , q^{\frac{1}{2}}
    z_{m-2},
    q \, z_{m-1})
    \\
    + \dots
    +
    q^{2(m-a-2)} \,  z_{a+2}^{~2} \cdots z_{m-1}^{~2} \,
    H_m^{(a)}(q \, x, y, z_1, \dots, z_a,q^{\frac{1}{2}}z_{a+1},
    q \, z_{a+2}, \dots, q \, z_{m-1})
    \\
    +
    q^{2(m-a-1)} \, z_{a+1}^{~2} \cdots z_{m-1}^{~2} \,
    \Bigl(
    H_m^{(a)}(q \, x, y, z_1,\dots,z_{a-1},q^{\frac{1}{2}} z_a,
    q \, z_{a+1} , \dots, q \, z_{m-1})
    \\
    +
    q \, z_a^{~2} \,
    H_m^{(a)}(q \, x, y, z_1,\dots,z_{a-2},q^{\frac{1}{2}} z_{a-1},
    q \, z_{a} , \dots, q \, z_{m-1})
    \\
    +
    \dots +
    q^a \, z_1^{~2} \cdots z_a^{~2} \,
    H_m^{(a)}(q \, x, y, q \, z_1,
    \dots, q \, z_{m-1})
    \Bigr)
    \label{gen_na_2}
  \end{multline}

  We set
  $z_1=z_2=\dots=z_{m-1}=z$ in
  eqs.~\eqref{gen_na_1}
  and~\eqref{gen_na_2}.
  Combining these two equations, we obtain
  \begin{multline}
    \label{dif_H_m_a_gen}
    -(1-y) \, H_m^{(a)}(x,y, \boldsymbol{z}) -
    (x \, y - q \, z^2) \,  H_m^{(a)}(x, q\,y, \boldsymbol{z})
    + 1- q^{a+1} \, z^{2a+2}
    \\
    =
    q^{2m-a} (1-x) \, y \, z^{2m} \,
    H_m^{(a)}(q \, x, q\, y, q\, \boldsymbol{z}) 
  \end{multline}
  When we set
  $(x,y,z) \to (q\, x, x, x)$, and by definition~\eqref{H_xyz_H}
  we can conclude that $H_m^{(a)}(x)$ satisfies
  eq.~\eqref{difference_general_a}.
\end{proof}

\begin{coro}
  Let the function $\Tilde{H}_m^{(a)}(x)$ be
  \begin{equation*}
    \Tilde{H}_m^{(a)}(x)
    =
    \frac{1}{x}\,
    \left(
      H_m^{(a)}( x, x, q^{-\frac{1}{2}} \boldsymbol{x})
      -
      1
    \right)
  \end{equation*}
  It satisfies the difference equation;
  \begin{equation}
    \Tilde{H}_m^{(a)}(x)
    =
    1+x+\dots + x^{2a}
    - q^{m-a} \, x^{2m}
    -
    q^{m-a+1} \, x^{2m+1} \,
    \Tilde{H}_m^{(a)}(q \, x)
  \end{equation}
\end{coro}

\begin{proof}
  Setting
  $(x,y,z) \to (x, x, q^{-\frac{1}{2}} x)$ in
  eq.~\eqref{dif_H_m_a_gen} gives above equation.
\end{proof}


We see from the right hand side of eq.~\eqref{H_m_a_and} that
\begin{equation}
  H_m^{(a)}(x=1)
  =
  (q^{a+1}, q^{2m-a}, q^{2m+1} ; q^{2m+1})_\infty
  \label{H_m_a_unity}
\end{equation}
For the left hand side we have the following identity which follows
from eq.~\eqref{binomial_3};
\begin{lemma}
\begin{multline}
  \label{H_m_a_finite}
  H_m^{(a)}(x)
  =
  (q \, x)_\infty
  \sum_{k_1,k_2,\dots,k_{m-1}=0}^\infty
  \frac{
    q^{k_1^{~2} + \dots + k_{m-1}^{~2} + k_{a+1} + \dots + k_{m-1}}
  }{
    (x \, q)_{k_{m-1}}
  } \,
  x^{2 \sum_{i=1}^{m-1} k_i + k_{m-1}}
  \\
  \times
  \left(
    \prod_{\substack{
        i=1 \\
        i \neq a
      }}^{m-2}
    \begin{bmatrix}
      k_{i+1} \\
      k_i
    \end{bmatrix}
  \right) \,
  \begin{bmatrix}
    k_{a+1} + 1 \\
    k_a
  \end{bmatrix}
  \\
  + (1 - x) \,
  \sum_{k_1,k_2,\dots,k_m=0}^\infty
  \bigl(
  (q \, x)_{k_m} - (q \, x)_\infty
  \bigr) \, x^{k_m} \,
  \left(
    \prod_{i=1}^{a-1}
    q^{k_i^{~2}} \, x^{2 k_i} \,
    \begin{bmatrix}
      k_{i+1} \\
      k_i
    \end{bmatrix}
  \right)
  \\
  \times
  q^{k_a^{~2}} \, x^{2 k_a} \,
  \begin{bmatrix}
    k_{a+1} + 1 \\
    k_a
  \end{bmatrix}
  \cdot
  \left(
    \prod_{i=a+1}^{m-1}
    q^{k_i^{~2} + k_i} \, x^{2 k_i} \,
    \begin{bmatrix}
      k_{i+1} \\
      k_i
    \end{bmatrix}
  \right)
\end{multline}
\end{lemma}

To check eq.~\eqref{H_m_a_unity}  from  eq.~\eqref{H_m_a_finite}, we
need a  variant of the Andrews--Gordon identity;
\begin{prop}
  \label{prop:A_and_Bailey}
  \begin{align}
    &
    \frac{
      (q^{a+1}, q^{2m-a}, q^{2m+1}; q^{2m+1})_\infty
    }{
      (q)_\infty
    }
    \nonumber 
    \\
    & =
    \sum_{k_1,\dots,k_{m-1}=0}^\infty
    \frac{
      q^{k_1^{~2} + \dots + k_{m-1}^{~2} + k_{a+1} + \dots + k_{m-1}}
    }{
      (q)_{k_{m-1}}
    }
    \,
    \left(
      \prod_{\substack{
          i=1 \\
          i \neq a
        }}^{m-2}
      \begin{bmatrix}
        k_{i+1} \\
        k_i
      \end{bmatrix}
    \right) \,
    \begin{bmatrix}
      k_{a+1} + 1 \\
      k_a
    \end{bmatrix}
    \nonumber
  \end{align}
\end{prop}

To prove this proposition
we need a certain limit of the Bailey lemma
(see, \emph{e.g.}, Ref.~\cite{Andre86,PPaul87a,GAndre01a});
\begin{prop}[Bailey lemma]
  \label{prop:Bailey}
  If for $n\geq 0$
  \begin{equation}
    \beta_n = \sum_{r=0}^n
    \frac{\alpha_r}{(q)_{n-r} \, (x \, q)_{n+r}}
  \end{equation}
  then
  \begin{equation}
    \beta_n^\prime
    =
    \sum_{r=0}^n
    \frac{\alpha_r^\prime}{(q)_{n-r} \, (x \, q)_{n+r}}
  \end{equation}
  where
  \begin{align}
    \alpha_r^\prime
    & =
    \frac{
      (\rho_1)_r \, (\rho_2)_r
    }{
      (x \,q / \rho_1)_r \, (x \, q / \rho_2)_r
    } \,
    \left(
      \frac{ x \,q}{\rho_1 \, \rho_2}
    \right)^r \, \alpha_r
    \\[2mm]
    \beta_n^\prime
    & =
    \sum_{j=0}^\infty
    \frac{
      (\rho_1)_j \, (\rho_2)_j \,
      \left(
        \frac{x \, q}{\rho_1 \, \rho_2}
      \right)_{n-j}
    }{
      (q)_{n-j} \,
      (x \,q / \rho_1)_n \,
      (x \,q / \rho_2)_n
    } \,
    \left(
      \frac{x \, q}{\rho_1 \, \rho_2}
    \right)^j \,
    \beta_j
  \end{align}
\end{prop}

\begin{coro}
  \label{coro:Bailey}
  \begin{equation}
    \sum_{k=0}^n
    \frac{
      a_k \, x^k
    }{
      (q)_{n-k} \, (x \, q)_{n+k}
    }
    =
    \sum_{j=0}^n
    \frac{q^{j^2} \, x^j}{(q)_{n-j}}
    \sum_{k=0}^j
    \frac{
      a_k \, q^{-k^2}
    }{
      (q)_{j-k} \, (x \, q)_{j+k}
    }
  \end{equation}
\end{coro}  

\begin{proof}
We  take a limit $n\to \infty$ in Prop.~\ref{prop:Bailey}.
Then, we set $\rho_1=q^{-m}$, $\rho_2=q^{-n}$, and
$\alpha_k=q^{-k^2} a_k$
and take $m\to \infty$.
\end{proof}

\begin{coro}
  \begin{gather}
    \label{Bailey_1}
    \sum_{k= -\infty}^\infty
    \frac{
      c_k
    }{
      (q)_{n-k} \, (q)_{n+k}
    }
    =
    \sum_{j=0}^\infty
    \frac{q^{j^2}}{(q)_{n-j}}
    \sum_{k=-\infty}^\infty
    \frac{c_k \, q^{-k^2}}
    {(q)_{j-k} \, (q)_{j+k}}
    \\[2mm]
    \label{Bailey_2}
    \sum_{k= -\infty}^\infty
    \frac{
      c_k
    }{
      (q)_{n-k} \, (q)_{n+k-1}
    }
    =
    \sum_{j=0}^\infty
    \frac{q^{j^2-j}}{(q)_{n-j}}
    \sum_{k=-\infty}^\infty
    \frac{c_k \, q^{-(k^2-k)}}
    {(q)_{j-k} \, (q)_{j+k-1}}
  \end{gather}
\end{coro}

\begin{proof}
We set $x=1, q^{-1}$ in Corollary~\ref{coro:Bailey}
and take a symmetrization for $a_k$.
\end{proof}

\begin{proof}[Proof of Proposition~\ref{prop:A_and_Bailey}]
  We set the left hand side of Prop.~\ref{prop:A_and_Bailey} as
  $A_m^{(a)}$.
  We apply the triple Jacobi identity to $A_m^{(a)}$,
  and then use the Bailey chain
  recursively;
  \begin{align}
    A_m^{(a)}
    &
    \stackrel{\eqref{Jacobi}}{=}
    \frac{1}{(q)_\infty} \,
    \sum_{k= - \infty}^\infty
    (-1)^k \,
    q^{(m+\frac{1}{2}) k^2 - (m-a-\frac{1}{2}) k}
    \nonumber
    \\
    &
    \stackrel{\eqref{Bailey_2}}{=}
    \sum_{k_{m-1} \geq k_{m-2} \geq \dots \geq k_{a+1} \geq 0}^\infty
    \frac{
      q^{
        \sum_{i=a+1}^{m-1} ( k_i^{~2} - k_i)
      }
    }{
      (q)_{k_{m-1}-k_{m-2}} \cdots (q)_{k_{a+2} - k_{a+1}}
    }
    \nonumber 
    \\
    & \qquad \qquad \qquad \times
    \sum_{k= -\infty}^\infty
    (-1)^k \,
    \frac{q^{(a+\frac{3}{2}) k^2 - \frac{1}{2} k}}
    {
      (q)_{k_{a+1} - k}  \, (q)_{k_{a+1}+k-1}
    }
  \end{align}
  Here we have for arbitrary $c$ that
  \begin{multline*}
    \sum_{k=-\infty}^\infty
    \begin{bmatrix}
      2 \, n-1 \\
      n-k
    \end{bmatrix}
    \, (-1)^k \,
    q^{c k^2  - \frac{1}{2} k}
    \\
    \stackrel{\eqref{binomial_2}}{=}
    \sum_{k=-\infty}^\infty
    \begin{bmatrix}
      2 \, n \\
      n-k
    \end{bmatrix}
    \, (-1)^k \,
    q^{c k^2 - \frac{1}{2} k}
    -
    q^n
    \sum_{k=-\infty}^\infty
    \begin{bmatrix}
      2 \, n  -1 \\
      n + k
    \end{bmatrix}
    \, (-1)^k \,
    q^{c k^2 + \frac{1}{2} k}
  \end{multline*}
  which gives
  \begin{equation}
    \label{mid_relate}
    \sum_{k=-\infty}^\infty
    (-1)^k \,
    \frac{
      q^{c k^2 - \frac{1}{2} k}
    }{
      (q)_{n-k} \, (q)_{n+k-1}
    }
    =
    (1-q^n)
    \sum_{k=-\infty}^\infty
    (-1)^k \,
    \frac{
      q^{ c k^2 - \frac{1}{2} k}
    }{
      (q)_{n-k} \, (q)_{n+k}
    }
  \end{equation}
  Then we have
  \begin{align*}
    A_m^{(a)}
    &
    \stackrel{\eqref{mid_relate}}{=}
    \sum_{k_{m-1} \geq \dots \geq k_{a+1} \geq 0}^\infty
    \frac{
      q^{
        \sum_{i=a+1}^{m-1} ( k_i^{~2} - k_i)
      }
    }{
      (q)_{k_{m-1}-k_{m-2}} \cdots (q)_{k_{a+2} - k_{a+1}}
    }    \,
    (1-q^{k_{a+1}})
    \\
    & \qquad \qquad \qquad \times
    \sum_{k= -\infty}^\infty
    (-1)^k \,
    \frac{q^{(a+\frac{3}{2}) k^2 - \frac{1}{2} k}}
    {
      (q)_{k_{a+1} - k} \, (q)_{k_{a+1}+k}
    }
    \\
    &
    \stackrel{\eqref{Bailey_1}}{=}
    \sum_{k_{m-1}\geq \dots \geq k_1 \geq 0}^\infty
    \frac{
      q^{
        \sum_{i=a+1}^{m-1} ( k_i^{~2} - k_i)
        +
        \sum_{k=1}^a k_i^{~2}
      }
    }{
      (q)_{k_{m-1}-k_{m-2}} \cdots (q)_{k_2 - k_1} \,
      (q)_{k_1}
    } \,
    (1 - q^{k_{a+1}})
  \end{align*}
  We note that,
  in the last equality, we have also used
  \begin{equation*}
    \sum_{k=-\infty}^\infty
    \frac{
      (-1)^k \, q^{\frac{1}{2} k^2 - \frac{1}{2} k}
    }{
      (q)_{n-k} \, (q)_{n+k}
    }
    =
    \delta_{n,0}
  \end{equation*}
  After we  shift parameters;
  $(k_{a+1},\dots,k_{m-1}) \to
  (k_{a+1} +1, \dots, k_{m-1}+1)$,
  we get a statement of    Prop.~\ref{prop:A_and_Bailey}.
\end{proof}

\begin{proof}[Proof of Theorem~\ref{thm:general_m}]
We differentiate eqs.~\eqref{H_m_a_finite} and~\eqref{H_m_a_and} with
respect to $x$, and substitute $x\to 1$.
We obtain
\begin{multline}
  \label{identity_X}
  (q^{a+1}, q^{2m-a} , q^{2m+1} ; q^{2m+1})_\infty \,
  \left(
    \frac{2 \, m - 2 \, a - 1}{2}
    -
    \sum_{i=1}^\infty
    \frac{q^i}{1-q^i}
  \right)
  \\
  +
  (q)_\infty
  \sum_{k_{1}, \dots, k_{m-1}=0}^\infty
  q^{k_1^{~2} + \dots + k_{m-1}^{~2} + k_{a+1}+ \dots + k_{m-1}}
  \,
  \left(
    2 \sum_{i=1}^{m-1} k_i + k_{m-1}
    +
    \sum_{i=1}^{k_{m-1}}
    \frac{q^i}{1-q^i}
  \right)
  \\
  \times
  \frac{1-q^{k_{a+1}+1}}
  {
    (q)_{k_{m-1}-k_{m-2}} \cdots (q)_{k_{a+2}-k_{a+1}} \,
    (q)_{k_{a+1}-k_a+1} \, (q)_{k_a - k_{a-1}} \cdots
    (q)_{k_1}
  }
  \\
  -
  \sum_{k_{1}, \dots, k_{m-1}=0}^\infty
  \bigl(
  (q)_{k_m} - (q)_\infty
  \bigr) \,
  q^{k_1^{~2} + \dots + k_{m-1}^{~2} + k_{a+1}+ \dots + k_{m-1}}
  \left(
    \prod_{\substack{i=1 \\ i \neq a}}^{m-1}
    \begin{bmatrix}
      k_{i+1} \\
      k_i
    \end{bmatrix}
  \right) \cdot
  \begin{bmatrix}
    k_{a+1} + 1 \\
    k_a
  \end{bmatrix}
  \\
  =
  \sum_{n=0}^\infty
  \frac{n}{2} \,
  \chi_{8m+4}^{(a)}(n) \,
  q^{
    \frac{n^2 - (2m -2a-1)^2}{8 (2m+1)}
  }
\end{multline}
We substitute $q\to \mathrm{e}^{-t}$, and find eq.~\eqref{half_AG}.
As a result we obtain Theorem~\ref{thm:general_m}
applying
Prop.~\ref{prop:Mellin}.
\end{proof}

\section{Knot Invariant and Nearly Modular Form}
\label{sec:knot}

We comment on a relationship between our  $q$-series and the knot
invariant.
Generally the $q$-series $X_m^{(a)}(q)$ does not converge in any open
set of $q$, but it reduces to the finite number in a case of $q$ being
root of unity.
Furthermore this finite value coincides with Kashaev's invariant
(or, the colored Jones polynomial with a specific value) for the
($2\, m+1,2$)-torus knot.

We prepare the following $q$-series identity;
\begin{lemma}
  We set $a \geq c \geq 0$. Then we have
  \begin{equation}
    \label{bc_lemma}
    \sum_{b=c}^a
    (-1)^{b+c} \, 
    \frac{q^{\frac{1}{2}b(b+1)+ \frac{1}{2}c(c+1)}}
    {(q)_{a-b} \, (q)_{b-c}}
    =
    q^{c^2+c} \,
    \begin{bmatrix}
      a \\
      c
    \end{bmatrix}
  \end{equation}
\end{lemma}
\begin{proof}
  Using $b=c+n$, we have
  \begin{align*}
    \text{l.h.s.}
    &=
    q^{c^2+c} \sum_{n=0}^{a-c}
    (-1)^n \, 
    \frac{q^{\frac{1}{2}(n^2+n) + c n}}{(q)_{a-c-n} \,  (q)_n}
    \\
    & =
    \frac{q^{c^2+c}}{(q)_{a-c}}
    \sum_{n=0}^{a-c}
    \frac{(q^{-(a-c)})_n}
    {(q)_n} \,
    q^{n(1+a)}
    \\
    &
    \stackrel{\eqref{binomial_3}}{=}
    \frac{q^{c^2+c}}{(q)_{a-c}} \cdot
    \frac{(q^{c+1})_\infty}{(q^{a+1})_\infty}
    =
    q^{c^2+c} \,
    \frac{(q^{c+1})_{a-c}}{(q)_{a-c}}
  \end{align*}
  which proves eq.~\eqref{bc_lemma}.
\end{proof}

We set $N\in \mathbb{Z}_{>0}$ and define\footnote[4]{
We use $\mathrm{i}=\sqrt{-1}$.}
\begin{equation}
  \omega=\exp(2 \, \pi \, \mathrm{i}/N)
\end{equation}
Hereafter we mean that
\begin{gather*}
  (\omega)_n
  = \prod_{a=1}^n ( 1- \omega^a)
  \\[2mm]
  \begin{bmatrix}
    n \\
    c
  \end{bmatrix}
  =
  \frac{(\omega)_n}{(\omega)_c \, (\omega)_{n-c}}
\end{gather*}
and $^{*}$ denotes a complex conjugate.

\begin{prop}
  When  $q$ being the $N$-th root of unity,
  $X_m^{(0)}(q=\omega)$ defined in eq.~\eqref{define_X_a}
  coincides  with Kashaev's invariant for the
  ($2,2\,m+1$)-torus knot.
\end{prop}

\begin{proof}
  We  take an example for $m=2$ case;
  \allowdisplaybreaks
  \begin{align*}
    X_2^{(0)}(\omega)
    & =
    \sum_{0\leq k_1 \leq k_2 \leq N-1}
    (\omega)_{k_2} \,
    \omega^{k_1^{~2} + k_1} \,
    \begin{bmatrix}
      k_2 \\
      k_1
    \end{bmatrix}
    \\
    &
    \stackrel{\eqref{bc_lemma}}{=}
    \sum_{0 \leq k_1 \leq n \leq k_2 \leq N-1}
    (-1)^{n+k_1} \,
    (\omega)_{k_2} \,
    \frac{
      \omega^{\frac{1}{2} n (n+1) + \frac{1}{2} k_1 (k_1+1)}
    }{
      (\omega)_{k_2 - n} \, (\omega)_{n-k_1}
    }
    \\
    & =
    \sum_{0 \leq k_1 \leq n \leq k_2 \leq N-1}
    \frac{
      (\omega)_{k_2}
    }{
      (\omega)_{k_2 - n} \, (\omega)_{n-k_1}^{~*}
    } \,
    \omega^{k_1 ( 1+n)}
    \\
    &
    \stackrel{\eqref{omega_2}}{=}
    \sum_{0 \leq k_1 \leq n \leq N-1}
    \frac{N}{
      (\omega)_{n-k_1}^{~*}
    } \,
    \omega^{k_1 (1+n)}
    \\
    &
    \stackrel{\eqref{omega_1}}{=}
    \sum_{0 \leq k_1 \leq n \leq N-1}
    (\omega)_{N-n+k_1-1} \,
    \omega^{k_1 (1+n)}
    \\
    & =
    \sum_{\substack{
        a,b=0 \\
        0 \leq a+b \leq N-1
      }}^{N-1}
    (\omega)_{a+b} \, \omega^{-a b}
  \end{align*}
  where we have used parameters
  $a=k_1$,  $b=N-1-n$.
  The last expression coincides with Kashaev's invariant for the
  ($5,2$)-torus knot given in Ref.~\cite{KHikami02b}.
  Here we have used for $0 \leq a \leq N-1$
  \begin{gather}
    \label{omega_2}
    \sum_{a \leq b \leq N-1}
    \frac{(\omega)_b}{(\omega)_{b-a}} = N
    \\[2mm]
    \label{omega_1}
    (\omega)_{N-1} =
    (\omega)_a^{~*} \, (\omega)_{N-1-a}=
    N
  \end{gather}
  We can check for other  $m$'s
  that $X_m^{(0)}(\omega)$ reduces to  the invariant for the torus
  knot
  given in
  Ref.~\cite{KHikami02b}.
\end{proof}

As was proved in Ref.~\cite{KHikami02b}, the asymptotic expansion of
knot invariant in a limit $N\to \infty$ can be written explicitly;
  \begin{multline}
    X_m^{(0)}(\omega)
    \\
    \simeq
      \frac{2}{\sqrt{2 \, m+1}} \, N^{\frac{3}{2}} \,
      \mathrm{e}^{
        \frac{\pi \mathrm{i}}{4}
        -
        \frac{\pi \mathrm{i}}{N} \, \frac{(2m-1)^2}{4 (2m+1)}
      } \,
      \sum_{k=0}^{m-1} (-1)^k \, (m-k) \,
      \sin
      \left( \frac{2 \, k+1}{2 \, m+1} \, \pi \right) \,
      \mathrm{e}^{- N \pi \mathrm{i} \frac{(2 k+1)^2}{4 (2m+1)}
      }
      \\
      +
      \mathrm{e}^{
        - \frac{\pi \mathrm{i}}{N} \,
        \frac{(2m-1)^2}{4 (2m+1)}
      } \,
      \sum_{n=0}^\infty
      \frac{T_m^{(0)}(n)}{n!} \,
      \left(
        \frac{\pi}{4 \, (2 \, m+1) \, N \, \mathrm{i}}
      \right)^n  .
      \label{formula_2m}
  \end{multline}

A case of $m=1$ is given in Ref.~\cite{DZagie01a} as ``Kontsevich's
conjectural asymptotic formula''.
In proving above asymptotic expansion, we  used a previously known
another expression
for the  colored Jones invariant for the torus knot.
Correspondingly
for a  case of $a \neq 0$
we have a theorem;
\begin{theorem}[Conjecture in Ref.~\cite{KHikami02b}]
  \label{conj:main}
  \begin{multline}
    \label{asymptotic_X_omega}
    X_m^{(a)}(\omega)
    \\
    \simeq
    \frac{2}{\sqrt{2 \, m+1}} \,
    N^{\frac{3}{2}} \,
    \mathrm{e}^{
      \frac{\pi \mathrm{i}}{4}
      -
      \frac{\pi \mathrm{i}}{N}
      \frac{(2m-2 a -1)^2}{4 (2m+1)}
    } \,
    \sum_{k=0}^{m-1}
    (-1)^k \, (m-k) \,
    \sin \left( (a+1) \, \frac{2 \, k+1}{2 \, m+1} \, \pi \right)
    \,
    \mathrm{e}^{
      - N \pi \mathrm{i} \frac{(2k+1)^2}{4 (2m+1)}
    }
    \\
    +
    \mathrm{e}^{
      -
      \frac{\pi \mathrm{i}}{N}
      \frac{(2m-2 a -1)^2}{4 (2m+1)}
    } \,
    \sum_{n=0}^\infty
    \frac{
      T_m^{(a)}(n)}{n!}
    \,
    \left(
      \frac{\pi}{
        4 \, (2 \, m+1) \, N \, \mathrm{i}
      }
    \right)^n 
  \end{multline}
\end{theorem}

This theorem indicates a nearly modular property~\cite{DZagie01a}
of the function
$X_m^{(a)}(q)$  with weight $1/2$;
we define functions $\widetilde{\Phi}_m^{(a)}(\alpha)$ by
\begin{align}
  \widetilde{\Phi}_m^{(a)}(\alpha)
  &
  =
  \mathrm{e}^{
    \frac{(2m-2a-1)^2}{4 (2 m+1)}  \pi \mathrm{i} \alpha
  } \,
  X_m^{(a)}(\mathrm{e}^{2 \pi \mathrm{i} \alpha}) 
\end{align}
and
introduce a vector $\widetilde{\boldsymbol{\Phi}}_m(\alpha)$
by
\begin{equation}
  \widetilde{\boldsymbol{\Phi}}_m(\alpha)
  =
  \begin{pmatrix}
    \widetilde{\Phi}_m^{(m-1)}(\alpha) \\
    \vdots \\
    \widetilde{\Phi}_m^{(1)}(\alpha) \\
    \widetilde{\Phi}_m^{(0)}(\alpha) \\
  \end{pmatrix} 
\end{equation}
Eq.~\eqref{asymptotic_X_omega} is reformulated into
\begin{equation}
  \label{modular_Phi_tilde}
  \widetilde{\boldsymbol{\Phi}}_m
  \left(\frac{1}{N}\right)
  +
  (- \mathrm{i} \, N)^{\frac{3}{2}} \,
  \mathbf{M}^{(2m+1)} \,
  \widetilde{\boldsymbol{\Phi}}_m(-N)
  =
  \sum_{n=0}^\infty
  \frac{  \mathbf{T}_m(n) }{n!} \,
  \left(
    \frac{\pi}{4 \, (2m+1)  \, \mathrm{i} \, N}
  \right)^n 
\end{equation}
where $\mathbf{M}^{(2m+1)}$ is an $m\times m$ matrix with an entry
\begin{align}
  \label{define_matrix}
  \Bigl(\mathbf{M}^{(2m+1)}\Bigr)_{1\leq a, b \leq m}
  & =
  \frac{2}{\sqrt{2 \, m+1}} \,
  \cos
  \left(
    \frac{(2 \, a-1) \, (2 \, b-1)}{2 \, (2 \, m+1)} \, 
    \pi
  \right)
\end{align}
and
\begin{equation*}
  \mathbf{T}_m(n)
  =
  \begin{pmatrix}
    T_m^{(m-1)}(n)
    \\
    \vdots
    \\
    T_m^{(1)}(n)
    \\
    T_m^{(0)}(n)
  \end{pmatrix} 
\end{equation*}

\begin{proof}[Proof of Theorem~\ref{conj:main}]
  Eq.~\eqref{identity_X} indicates that  the function
  $\widetilde{\Phi}_m^{(a)}(\alpha)$ coincides with
  a limit
  $\tau \to \alpha \in \mathbb{Q}$
  of the $q$-series~\eqref{half_AG}
  \begin{equation}
    \widetilde{\Phi}_m^{(a)}(\tau)
    =
    -\frac{1}{2}
    \sum_{n=0}^\infty n \, \chi_{8 m+4}^{(a)}(n) \,
    q^{\frac{1}{8 (2m+1)} n^2}
  \end{equation}
  for   $q=\mathrm{e}^{2 \pi \mathrm{i} \tau}$.
  This is regarded as the Eichler integral of
  \begin{equation}
    \Phi_m^{(a)}(\tau)
    =
    \sum_{n=0}^\infty \chi_{8 m+4}^{(a)}(n) \,
    q^{\frac{1}{8 (2m+1)} n^2}
  \end{equation}
  which is modular with weight $1/2$;
  it is straightforward to see that
  \begin{gather}
    \Phi_m^{(a)}(\tau +1) =
    \mathrm{e}^{
      \frac{ (2m-2a-1)^2}{4(2m+1)} \pi \mathrm{i}
    } \,
    \Phi_m^{(a)}(\tau)
  \end{gather}
  and using the Poisson summation formula we obtain
  \begin{gather}
    \boldsymbol{\Phi}_m (\tau)
    =
    \sqrt{
      \frac{\   \mathrm{i}  \   }{\tau}
    } \,
    \mathbf{M}^{(2m+1)} \,
    \boldsymbol{\Phi}_m (-1/\tau)
  \end{gather}
  where $\mathbf{M}^{(2m+1)}$ is an $m\times m$ matrix defined in
  eq.~\eqref{define_matrix}, and
  \begin{equation*}
    {\boldsymbol{\Phi}}_m(\tau)
    =
    \begin{pmatrix}
      {\Phi}_m^{(m-1)}(\tau) \\
      \vdots \\
      {\Phi}_m^{(1)}(\tau) \\
      {\Phi}_m^{(0)}(\tau) \\
    \end{pmatrix} 
  \end{equation*}

  To prove eq.~\eqref{modular_Phi_tilde} following
  Refs.~\cite{LawrZagi99a,DZagie01a},
  we study
  an analogue of the Eichler integral defined by
  \begin{equation}
    \widehat{\Phi}_m^{(a)}(z)
    =
    \frac{\sqrt{(2 \, m+1) \, \mathrm{i}}}{2 \, \pi}
    \,
    \int_{z^*}^\infty
    \frac{\Phi_m^{(a)}(\tau)}{
      (\tau - z )^{\frac{3}{2}}
    } \,
    \mathrm{d} \tau
  \end{equation}
  which is defined for $z$ in the lower half plane $z \in \mathbb{H}^-$.
  By performing an integration term by term,
  we have
  \begin{align*}
    \widehat{\Phi}_m^{(a)}
    (z)
    & =
    \frac{\sqrt{(2 \, m+1) \, \mathrm{i}}}{2 \, \pi}
    \,
    \sum_{n=0}^\infty
    \chi_{8 m+4}^{(a)}(n) \,
    \int_{z^*}^\infty
    \frac{
      \mathrm{e}^{\pi \mathrm{i} \tau \frac{n^2}{4 (2m+1)}}
    }{
      (\tau - z)^{3/2}
    } \,
    \mathrm{d} \tau
    \\
    & \stackrel{z \to \alpha \in \mathbb{Q}}{\to}
    -\frac{1}{2}
    \sum_{n=0}^\infty 
    n \, \chi_{8m+4}^{(a)}(n) \,
    \mathrm{e}^{\frac{n^2}{4 (2m+1)} \pi \mathrm{i} \alpha}
  \end{align*}
  which shows
  \begin{equation}
    \widehat{\Phi}_m^{(a)}(\alpha)
    =
    \widetilde{\Phi}_m^{(a)}(\alpha)
    \label{hat_and_tilde}
  \end{equation}
  Note that l.h.s. is a limiting value
  from a lower half plane $\mathbb{H}^-$
  while r.h.s. is given from an upper half plane $\mathbb{H}$.

  To see a  modular property of
  $\widehat{\boldsymbol{\Phi}}_m(\tau)$,
  we define
  the period function by
  \begin{equation}
    r_m^{(a)}(z; \alpha)
    =
    \frac{\sqrt{(2 \, m+1) \, \mathrm{i}}}{2 \, \pi}
    \,
    \int_\alpha^\infty
    \frac{\Phi_m^{(a)}(\tau)}{
      (\tau - z )^{\frac{3}{2}}
    } \,
    \mathrm{d} \tau
  \end{equation}
  where $\alpha \in \mathbb{Q}$.
  It is defined for  $z \in \mathbb{H}^-$,
  but it is analytically continued to $\mathbb{R}$.
  We then have
  \begin{align}
    \sum_{b=1}^{m}
    \left( \mathbf{M}^{(2m+1)} \right)_{a,b}
    \,
    \widehat{\Phi}_m^{(m-b)}(-1/z)
    & =
    \sum_{b=1}^{m}
    \left( \mathbf{M}^{(2m+1)} \right)_{a,b}
    \,
    \frac{\sqrt{(2 \, m+1) \, \mathrm{i}}}{2 \, \pi}
    \,
    \int_{z^*}^0
    \frac{\Phi_m^{(m-b)}(-1/s)}{
      (-s^{-1} + z^{-1})^{3/2}
    }
    \frac{\mathrm{d} s}{s^2}
    \nonumber
    \\
    & =
    - ( \mathrm{i} \, z)^{3/2} \,
    \frac{\sqrt{(2 \, m+1) \, \mathrm{i}}}{2 \, \pi}
    \,
    \int_{z^*}^0
    \frac{\Phi_m^{(m-a)}(-1/s)}{
      (s -  z)^{3/2}
    }
    \mathrm{d} s
    \nonumber
    \\
    & =
    - ( \mathrm{i} \, z)^{3/2} \,
    \left(
      \widehat{\Phi}_m^{(m-a)}(z)
      -
      r_m^{(m-a)}(z;0)
    \right)
    \label{modular_hat}
  \end{align}
  We consider   a limit $z\to 1/N$ in eq.~\eqref{modular_hat}.
  We see that
  an asymptotic expansion of $r_m^{(a)}(1/N;0)$ in $N\to\infty$
  is given by
  \begin{align*}
    r_m^{(a)}\left(\frac{1}{N}; 0 \right)
    & =
    - \frac{1}{2}
    \sum_{k=0}^\infty
    \frac{
      L(-2 \, k -1 , \chi_{8m+4}^{(a)})
    }{k!}
    \,
    \left(
      \frac{ \pi \, \mathrm{i}}{4 \, (2 \, m+1) \, N}
    \right)^k
  \end{align*}
  and 
  eq.~\eqref{hat_and_tilde} indicates that
  $\widehat{\Phi}_m^{(a)}(z)$  coincides with
  $\widetilde{\Phi}_m^{(a)}(z)$ for $z=-N$ and $z=1/N$.
  Recalling eq.~\eqref{T_and_L_function},
  we can conclude eq.~\eqref{modular_Phi_tilde}.
\end{proof}

\end{document}